\documentclass[12pt]{article}

\usepackage{slashbox}
\usepackage{amssymb}
\usepackage{amsmath}
\usepackage{amscd}
\usepackage[dvips]{graphicx}

\numberwithin{equation}{section}

\newtheorem{theorem}{Theorem}[section]

\begin{document}

\title{ASYMPTOTIC EXPANSIONS FOR THE CONDITIONAL SOJOURN TIME DISTRIBUTION IN THE $M/M/1$-PS QUEUE}

\author{
Qiang Zhen\thanks{
Department of Mathematics, Statistics, and Computer Science,
University of Illinois at Chicago, 851 South Morgan (M/C 249),
Chicago, IL 60607-7045, USA.
{\em Email:} qzhen2@uic.edu.}
\and
and
\and
Charles Knessl\thanks{
Department of Mathematics, Statistics, and Computer Science,
University of Illinois at Chicago, 851 South Morgan (M/C 249),
Chicago, IL 60607-7045, USA.
{\em Email:} knessl@uic.edu.\
\newline\indent\indent{\bf Acknowledgement:} This work was partly supported by NSF grant DMS 05-03745.
}}
\date{June 27, 2007 }
\maketitle

\begin{abstract}
\noindent  We consider the $M/M/1$ queue with processor sharing. We study the conditional sojourn time distribution, conditioned on the customer's service requirement, in various asymptotic limits. These include large time and/or large service request, and heavy traffic, where the arrival rate is only slightly less than the service rate. The asymptotic formulas relate to, and extend, some results of Morrison \cite{MO} and Flatto \cite{FL}.
\end{abstract}

\newpage

\section{Introduction}
One of the most interesting service disciplines within queueing theory is processor sharing (PS). Unlike other disciplines such as first-in first-out (FIFO), under the PS discipline every customer gets a share of the server. The advantage of this discipline over FIFO is that customers requiring only short amounts of service get through the system more rapidly than with other service disciplines. The PS discipline was apparently introduced by Kleinrock \cite{KLanalysis}, \cite{KLtime}, and has been the subject of much further investigation over the past forty years.

In  \cite{COF}, Coffman, Muntz, and Trotter derived an expression for the Laplace transform of the distribution of the waiting time in the $M/M/1$-PS model. We shall denote the waiting time (in the steady state) by $\mathbf{W}$ and the sojourn time by $\mathbf{V}$, which is the total time from when a customer arrives to when that customer leaves the system, after being served. In the $M/M/1$-PS model we denote the Poisson arrival rate by $\lambda$ and the exponential service time density by $\mu e^{-\mu x}$. The traffic intensity is $\rho=\lambda/\mu$. If we are given a service time $x$, then we let $\mathbf{V}(x)$ be the conditional sojourn time. The ``waiting time" $\mathbf{W}(x)$ is defined by $\mathbf{W}(x)=\mathbf{V}(x)-x$, and with this notation $\mathbf{E}[e^{-s\mathbf{W}(x)}]$ is the Laplace transform of the conditional waiting time distribution derived in  \cite{COF}.

Using the results in  \cite{COF}, Morrison \cite{MO} studied the unconditional sojourn time distribution in the $M/M/1$-PS model, in the heavy traffic limit where $\lambda\uparrow\mu$ (thus $\rho\uparrow 1$). Let us denote the unconditional sojourn time density by $p(t)dt=\Pr\big[\mathbf{V}\in(t,t+dt)\big]$ and the conditional density by
\begin{equation}\label{eqI1}
p(t|x)dt=\Pr\big[\mathbf{V}(x)\in(t,t+dt)\big].
\end{equation}
For the $M/M/1$-PS model we have
\begin{equation}\label{eqI2}
p(t)=\int_0^t \mu e^{-\mu x}p(t|x)dx
\end{equation}
and we note that $p(t|x)$ has support for $t>x$, with a probability mass along $t=x$.

In  \cite{MO} Morrison assumes that $\lambda\uparrow\mu$ and defines $\epsilon=1-\lambda/\mu=1-\rho$. He obtains approximations to $p(t)$ on the time scales $t=O(1)$, $t=O(\epsilon^{-1})$, $t=O(\epsilon^{-2})$ and $t=O(\epsilon^{-3})$. Most of the mass, in this asymptotic limit, is concentrated in the range $t=O(\epsilon^{-1})$.

A service discipline seemingly unrelated to PS is random order service (ROS), where customers are chosen for service at random. The $M/M/1$-ROS model has been studied by many authors, see Vaulot \cite{VA}, Riordan \cite{RI}, Kingman \cite{KI} and Flatto \cite{FL}. In  \cite{FL} an explicit integral representation is derived for the waiting time distribution, from which the following tail behavior is computed:
\begin{equation}\label{eqI3}
\Pr\left[\mathbf{W}_{\mathrm{ROS}}>t\right]\sim e^{-At}e^{-Bt^{1/3}}Ct^{-5/6},\textrm{ }\textrm{ }\textrm{ }  t\rightarrow\infty.
\end{equation}
Here $A=(1-\sqrt{\rho})^2$ if we scale time to make the service rate $\mu=1$, and $\rho<1$. This formula appeared previously in the book of Riordan \cite{RIbook} (pg.105) and was apparently first obtained by Pollaczek \cite{PO}. Cohen \cite{COH} established a relationship between the sojourn time in the PS model and the waiting time in the ROS model,
\begin{equation}\label{eqI4}
\Pr\left[\mathbf{V}>t\right]=\frac{1}{\rho}\Pr[\mathbf{W}_{\mathrm{ROS}}>t],
\end{equation}
which extends also to the more general $G/M/1$ case. In  \cite{BO} relations of the form (\ref{eqI4}) are explored for other models, such as finite capacity queues, repairman problems and networks. Later in this paper we discuss the relationship of the tail formula in (\ref{eqI3}) to the heavy traffic expansions in  \cite{MO}.

In this paper we study the conditional sojourn time distribution for the $M/M/1$-PS model in two cases. First we consider a fixed $\rho<1$ and obtain expansions of $p(t|x)$ for $t$ and/or $x\rightarrow\infty$. From these (\ref{eqI3}) is readily obtained by using (\ref{eqI2}) for $t$ large. Then we consider the heavy traffic limit where $\rho\uparrow 1$, and again obtain approximations for several ranges of the space-time plane. From these results all of the expansion in  \cite{MO} can be recovered.

We mention some other work on PS queues. The $G/M/1$-PS model was studied by Ramaswami \cite{RA} and the $M/G/1$-PS model by Yashkov \cite{YAproc}, \cite{YAmath} and by Ott \cite{OT}. For the latter there is a complicated expression for the Laplace transform of the conditional sojourn time distribution. For the special case of $M/D/1$-PS this expression simplifies considerably, and then the tail behavior was obtained by Egorova, et. al. \cite{EG} in the form $\Pr[\mathbf{V}>t]\sim B^\prime e^{-A^\prime t}$. Comparing this to (\ref{eqI4}) with (\ref{eqI3}) the interesting question arises as to what are the variety of possible tail behaviors for the general $M/G/1$-PS model. For the $G/G/1$-PS queue the tail was shown by Mandjes and Zwart \cite{MA} to be roughly exponential, in that $\log\{\Pr[\mathbf{V}>t]\}\sim -A_0 t$ as $t\rightarrow\infty$. This assumes that the arrival and service distributions have exponentially small tails, so that their moment generating functions are analytic in some part of the left half-plane.

There have also been several recent investigations into tail behaviors of PS models that are ``heavy tailed", in that the service time distribution has an algebraic tail. The $M/G/1$ model with this assumption is analyzed by Zwart and Boxma \cite{ZW}, where it is shown that the sojourn time has a similar algebraic tail. A thorough recent survey of sojourn time asymptotics in PS models can be found in  \cite{BONU}, where both heavy and light tailed distributions are discussed, as well as various approaches to obtaining the asymptotics.

The remainder of the paper is organized as follows. In section 2 we summarize the main results, as Theorems 1 and 2. In section 3 and 4 we sketch the main points of the derivations. A brief discussion appears in section 5.

\section{Summary of results}
We denote the mean rate for the Poisson arrival process by $\lambda$, and assume that the mean required service time is unity. Then the traffic intensity is $\rho=\lambda/\mu=\lambda>0$, where $\mu=1$ as in the notation of Coffman, et.al. \cite{COF}.

It was shown in  \cite{COF}, under the stability condition $\rho<1$, that the Laplace transform of the equilibrium waiting time distribution, conditioned on the job requiring $x$ units of service, is given by
\begin{equation}\label{eqS1}
\mathbf{E}[e^{-s \mathbf{W}(x)}]=\frac{(1-\rho)(1-\rho r^2)e^{-\rho (1-r)x}}{(1-\rho r)^2-\rho (1-r)^2 e^{-x (1-\rho r^2)/r}},
\end{equation}
where
\begin{equation}\label{eqS2}
r=r(s;\rho)=\frac{1}{2\rho}\Big[(1+\rho+s)-\sqrt{(1+\rho+s)^2-4\rho}\Big].
\end{equation}

\noindent Taking the inverse Laplace transform of (\ref{eqS1}), the probability density of the sojourn time, conditioned on service time $x$, is
\begin{eqnarray}\label{eqS3}
p(t|x) &=& \frac{1}{2\pi i}\int_{Br} e^{s(t-x)}\mathbf{E}[e^{-s\mathbf{W}(x)}] ds\nonumber\\
        &=& \frac{1}{2\pi i} \int_{Br} \frac{e^{s t}e^{-(s+\rho-\rho r)x}(1-\rho)(1-\rho r^2)} {(1-\rho r)^2-\rho(1-r)^2 e^{-x(1-\rho r^2)/r}} ds.
\end{eqnarray}

\noindent Here $Br$ is a vertical contour in the complex $s$-plane, on which $\Re(s)>0$.
Note that this density contains the term $\delta(t-x)(1-\rho)e^{-\rho x}$, which corresponds to probability mass at $t=x$. The sojourn time is equal to the service time if the customer arrives to an empty system and there are no arrivals during that customer's service time.
Analyzing the integral (\ref{eqS3}), we obtain the following results:

\begin{theorem} \label{th1}
For $\rho<1$, the conditional sojourn time density has the following asymptotic expansions:
\begin{enumerate}

\item $x\rightarrow \infty$, $t-x\rightarrow 0^+$ with $x(t-x)=O(1)$,
\begin{equation}\label{eqth11}
p(t|x)\sim(1-\rho)\sqrt{\frac{\rho x}{t-x}}e^{-\rho x}I_1(2\sqrt{\rho x(t-x)})
\end{equation}
where $I_1(\cdot)$ is the modified Bessel function.

\item $x$,$t\rightarrow \infty$, $1<t/x<\infty$,
\begin{equation}\label{eqth12}
p(t|x)\sim \frac{(1-\rho)\rho^{1/4}e^{2t\sqrt{\rho(1-x/t)}}e^{-(1+\rho)t}e^x}{2\sqrt{\pi}\sqrt{x}\left[(t/x)(t/x-1)\right]^{3/4}\left[\sqrt{t/x}-\sqrt{\rho(t/x-1)}\right]^2}.
\end{equation}

\item $x\rightarrow\infty$, $t/x^2 = a = O(1)$,
\begin{eqnarray}\label{eqth136}
p(t|x)\sim&& \frac{(1+\sqrt{\rho})e^{-(1-\sqrt{\rho})^2t}e^{(1-\sqrt{\rho})x}}{2\sqrt{\pi}\rho^{1/4}(1-\sqrt{\rho})a^{5/2}x^3}\nonumber\\
          &&\cdot\sum_{n=0}^{\infty}e^{-(2n+1)^2\sqrt{\rho}/(4a)}\left[(2n+1)^2\sqrt{\rho}-2a\right].
\end{eqnarray}
By the Poisson summation formula, we also have
\begin{equation}\label{eqth137}
p(t|x)\sim \frac{2\pi^2(1+\sqrt{\rho})e^{-(1-\sqrt{\rho})^2t}e^{(1-\sqrt{\rho})x}}{\rho(1-\sqrt{\rho})x^3}\sum_{m=1}^{\infty}(-1)^{m+1}m^2 e^{-\pi^2m^2a/\sqrt{\rho}}.
\end{equation}

\item $x=O(1)$, $t\rightarrow\infty$,
\begin{equation}\label{eqth148}
p(t|x)\sim F(x) e^{r_*(x)t},
\end{equation}
where for $x>x_*=\frac{1-\sqrt{\rho}}{\sqrt{\rho}(1+\sqrt{\rho})}$,
\begin{eqnarray*}
F(x) =&& \frac{2\sqrt{\rho}(1-\rho) e^{(1-\sqrt{\rho}\cos v)x}\sin^2 v\cos(\sqrt{\rho} x\sin v)}{(1+\rho-2\sqrt{\rho}\cos v)^2\sqrt{\rho} x \cos v+(1-\rho)(1+\rho-2\sqrt{\rho}\cos v)}\\
     &\cdot& \left[(1-\rho)\sin v\tan(\sqrt{\rho}x\sin v)+2\sqrt{\rho}-(1+\rho)\cos v\right]
\end{eqnarray*}
and
\begin{equation}\label{eqth149}
r_*(x)=-1-\rho+2\sqrt{\rho}\cos(v(x))
\end{equation}
where $v=v(x)$ is the smallest positive root of
\begin{equation}\label{eqth1410}
\left(\frac{1-\sqrt{\rho}e^{-iv}}{\sqrt{\rho}-e^{-iv}}\right)^2=e^{-2i\sqrt{\rho}x\sin v}.
\end{equation}
For $0<x<x_*$ we have
$$ r_*(x)=-1-\rho-2\sqrt{\rho}\cosh(u(x)) $$
where $u(x)$ is the unique positive root of 
$$ \left(\frac{1+\sqrt{\rho}e^{-u}}{\sqrt{\rho}+e^{-u}}\right)^2=e^{2\sqrt{\rho}x \sinh u} $$ 
and in $F(x)$ we must replace $(\sin v, \cos v)$ by $(i\sinh u,-\cosh u)$. When $x=x_*$ we have 
$r_*(x_*)=-(1+\sqrt{\rho})^2$, $v(x_*)=\pi$, $u(x_*)=0$ and 
$$F(x_*)=\frac{6\sqrt{\rho}(1+\sqrt{\rho})^2}{1+4\sqrt{\rho}+\rho}\exp\left(\frac{1}{\sqrt{\rho}}-1\right).$$

$r_*(x)$ has the following expansions, as $x\rightarrow\infty$:
\begin{equation}\label{eqth1411}
r_*(x)= -(1-\sqrt{\rho})^2-\frac{\pi^2}{\sqrt{\rho}x^2}+\frac{2\pi^2(1+\sqrt{\rho})}{\rho(1-\sqrt{\rho})x^3}+O(x^{-4}),
\end{equation}

and as $x\rightarrow 0^+$:
\begin{eqnarray}\label{eqth1412}
r_*(x)&=& \frac{\log\rho}{x}-\left[1+\rho+\frac{2(1-\rho)}{\log\rho}\right]\\
&&+\frac{2\rho(\log\rho)^2+(\rho^2-1)\log\rho-4(1-\rho)^2}{(\log\rho)^3}x+O(x^2)\nonumber.
\end{eqnarray}

\end{enumerate}
\end{theorem}

In the asymptotic matching region between cases $3$ and $4$ we have
\begin{eqnarray}\label{eqth1match}
p(t|x)\sim&& \frac{2\pi^2(1+\sqrt{\rho})}{\rho(1-\sqrt{\rho})x^3}\\
&&\cdot\exp\left\{(1-\sqrt{\rho})x+\left[-(1-\sqrt{\rho})^2-\frac{\pi^2}{\sqrt{\rho}x^2}+\frac{2\pi^2(1+\sqrt{\rho})}{\rho(1-\sqrt{\rho})x^3}\right]t\right\}\nonumber
\end{eqnarray}
which holds for $t\rightarrow\infty$ and $x\rightarrow\infty$, with $x=O(t^{1/3})$ or $t^{1/3}\ll x \ll t^{1/2}$.
By removing the condition on $x$ and noticing the relationship (\ref{eqI4}) between processor sharing and service in random order, we can recover Theorem 1.2 in Flatto \cite{FL}.
\\
\\
For the heavy-traffic case, where $\rho$ is close to 1, we let $\epsilon=1-\rho\rightarrow0^+$. Then we have the following results:

\begin{theorem} \label{th2}
For $\rho=1-\epsilon$, where $\epsilon\rightarrow0^+$, we let $t=T/\epsilon$ and $x=X/\epsilon$. The conditional sojourn time density has the following asymptotic expansions:

\begin{enumerate}
\item $x=O(1)$, $t=O(1)$,
\begin{equation}\label{eqth21}
p(t|x)\sim \frac{\epsilon}{2\pi i}\int_{Br} \frac{e^{st}e^{-(s+1-r_0)x}(1-r_0^2)}{(1-r_0)^2(1-e^{-x(1-r_0^2)/r_0})}ds,
\end{equation}
where
$$ r_0=\frac{1}{2}\left[s+2-\sqrt{s^2+4s}\right].$$

\item $x=O(1)$, $T=O(1)$,
\begin{equation}\label{eqth22}
p(t|x)=\frac{\epsilon}{x}e^{-T/x}-\frac{x+3}{6}\left[\delta(T)-\frac{2x-T}{x^2}e^{-T/x}\right]\epsilon^2+O(\epsilon^3).
\end{equation}

\item $X=O(1)$, $T=O(1)$, $T-X\rightarrow 0^+$ with $(T-X)=\epsilon^2 T_*=O(\epsilon^2)$,
\begin{equation}\label{eqth23}
p(t|x)\sim e^{-X/\epsilon}e^X\sqrt{X/T_*}I_1(2\sqrt{XT_*}).
\end{equation}

\item $X=O(1)$, $T=O(1)$ and $1<T/X<\infty$,
\begin{eqnarray}\label{eqth24}
p(t|x)\sim && \frac{\epsilon^{3/2}\left(\sqrt{T}+\sqrt{T-X}\right)^2}{2\sqrt{\pi}T^{3/4}(T-X)^{3/4}}\exp\left[T-\sqrt{T(T-X)}\right]\nonumber\\
&\cdot&\exp\left\{\frac{1}{\epsilon}\left[2\sqrt{T(T-X)}+X-2T\right]\right\}.
\end{eqnarray}

\item $X=\sqrt{\epsilon}Z=O(\sqrt{\epsilon})$, $T=O(1)$,
\begin{equation}\label{eqth25}
p(t|x)\sim \frac{2\epsilon^{3/2}}{\sqrt{\pi T}}\sum_{n=0}^\infty e^{-(2n+1)^2Z^2/(4T)},
\end{equation}
and by the Poisson summation formula,
\begin{equation}\label{eqth25psf}
p(t|x)\sim \frac{\epsilon^{3/2}}{Z}[1+2\sum_{n=1}^\infty (-1)^n e^{-n^2\pi^2 T/Z}].
\end{equation}

\item $X=O(1)$, $T=\Theta/\epsilon=O(\epsilon^{-1})$, we give the expansion in three different forms:
\begin{enumerate}
\item
\begin{equation}\label{eqth26a}
p(t|x)\sim \frac{4\epsilon^2e^{X/2}e^{-\Theta/4}}{\pi i}\int_{Br} \frac{e^{\xi\Theta}\sqrt{\xi}e^{-\sqrt{\xi}X}}{(1+2\sqrt{\xi})^2-(1-2\sqrt{\xi})^2e^{-2\sqrt{\xi}X}}d\xi.
\end{equation}

\item
\begin{eqnarray}\label{eqth26b}
&&p(t|x)\sim \frac{\epsilon^2}{\sqrt{2\pi}} \sum_{n=0}^\infty e^{(n+1)X}e^{-z_n^2/4}\sum_{l=0}^{2n}(-1)^l\frac{(2n)!}{l!(2n-l)!}(2\Theta)^{l/2}\nonumber\\
&&\cdot\left[\frac{4}{\sqrt{2\Theta}}D_{-l}(z_n)-4D_{-l-1}(z_n)+\sqrt{2\Theta}D_{-l-2}(z_n)\right].
\end{eqnarray}
Here $D_\nu(\cdot)$ is the parabolic cylinder function, and $z_n=\frac{(2n+1)X+\Theta}{\sqrt{2\Theta}}$.

\item
\begin{equation}\label{eqth26c}
p(t|x)\sim \epsilon^2\sum_{n=1}^{\infty}\widetilde{F}(X;v_n)e^{\widetilde{r}_*(X;v_n)\Theta}.
\end{equation}
Here
$$\widetilde{F}(X;v_n)=\frac{8v_n^2[(4v_n^2-1)\cos(Xv_n)+4v_n \sin(Xv_n)]}{(4v_n^2+1)[(4v_n^2+1)X+4]}e^{X/2},$$
$$\widetilde{r}_*(X;v_n)=-(v_n^2+1/4)-(v_n^2+1/4)\epsilon/2+O(\epsilon^2),$$

and $v_n=v_n(X)$ are the real positive roots of the equation
\begin{equation}\label{eqth26croot}
\left(\frac{2iv+1}{2iv-1}\right)^2=e^{-2iv X}.
\end{equation}

\end{enumerate}

If $\Theta=\sigma/\epsilon=O(\epsilon^{-1})$, we have
\begin{equation}\label{eqth26cmatch}
p(t|x)\sim \epsilon^2\widetilde{F}(X;v_1)e^{\widetilde{r}_*(X;v_1)\sigma/\epsilon}
\end{equation}
where $v_1=v_1(X)$ is the unique root in the interval $(0,\pi/X)$ of (\ref{eqth26croot}).

\end{enumerate}
\end{theorem}

By removing the condition on $x$ and using the results (\ref{eqth26cmatch}) and (\ref{eqth22}), we can recover the results in Morrison \cite{MO}.
The details are given in section 4. Here we also relate the unconditional tail expansion of Flatto \cite{FL} to the heavy traffic result of Morrison \cite{MO}.

\section{Brief derivations for the case $\rho<1$}
We assume that the traffic intensity $\rho$ is fixed and less than one. Consider the function (\ref{eqS1}) with (\ref{eqS2}). We first observe that if we replace $-\sqrt{(1+\rho+s)^2-4\rho}$ by $+\sqrt{(1+\rho+s)^2-4\rho}$ in (\ref{eqS2}), then $r$ becomes $1/(\rho r)$. But (\ref{eqS1}) is invariant under the transformation $r\rightarrow 1/(\rho r)$. From (\ref{eqS2}) and the fact that $r(s;\rho)$ appears in (\ref{eqS1}) it would seem that the integrand in (\ref{eqS3}) has branch points at $s=s_\pm$ with $s_*=-(1\pm\sqrt{\rho})^2=-1-\rho\mp2\sqrt{\rho}$, but from our discussion it follows that (\ref{eqS1}) has in fact no branch points. It has poles where the denominator in (\ref{eqS1}) vanishes. We also note that $r(s_-;\rho)=1/\sqrt{\rho}$ and $(1-\rho r^2)/r>0$ for $s$ real and $s>s_-$.

We define
\begin{equation}\label{eqD3_1}
\varphi(s)=\varphi\left(s;\frac{t}{x}\right)=s\left(\frac{t}{x}-1\right)+\rho r(s;\rho)
\end{equation}
and rewrite (\ref{eqS3}) as
\begin{eqnarray}\label{eqD3_2}
p(t|x)&=&\frac{(1-\rho)e^{-\rho x}}{2\pi i}\int_{Br}e^{x\varphi(s)}\frac{1-\rho r^2}{(1-\rho r)^2-\rho(1-r)^2e^{-x(1-\rho r^2)/r}}ds\nonumber\\
&\sim& \frac{(1-\rho)e^{-\rho x}}{2\pi i}\int_{Br}e^{x\varphi(s)}\frac{1-\rho r^2(s;\rho)}{\left[1-\rho r(s;\rho)\right]^2}ds
\end{eqnarray}
where the last approximation holds when $x\rightarrow\infty$ and $\Re(s)>s_-$ on the Bromwich contour $Br$, which is a vertical contour in the complex $s$-plane.

Consider first the limit $x,t\rightarrow\infty$ with $1<t/x<\infty$. The second integrand in (\ref{eqD3_2}) has saddle points where $\varphi^\prime(s)=0$ so that there is a saddle along the real axis at $s_0$, where
\begin{equation}\label{eqD3_3}
s_0=s_0\left(\frac{t}{x}\right)=-1-\rho+\frac{\sqrt{\rho}(2t-x)}{\sqrt{t(t-x)}}.
\end{equation}
Note that $s_-<s_0<\infty$ whenever $1<t/x<\infty$, and $s_0\rightarrow\infty$ as $t/x\downarrow 1$, while $s_0\rightarrow s_-$ as $t/x\rightarrow\infty$. If we shift $Br$ in (\ref{eqD3_2}) to $Br^\prime$, on which $\Re(s)=s_0$, and use the standard Laplace method (see, e.g., Wong \cite{WO}) we get
\begin{equation}\label{eqD3_4}
p(t|x)\sim\frac{(1-\rho)e^{-\rho x}\left[1-\rho r^2(s_0;\rho)\right]e^{x\varphi(s_0)}}{\sqrt{2\pi}\left[1-\rho r(s_0;\rho)\right]^2\sqrt{x\varphi^{\prime\prime}(s_0)}}.
\end{equation}
But from (\ref{eqD3_1}) and (\ref{eqS2}) we get
\begin{equation}\label{eqD3_5}
r(s_0;\rho)=\sqrt{\frac{t-x}{\rho t}},
\end{equation}
\begin{equation}\label{eqD3_6}
x\varphi(s_0)=(x-t)(1+\rho)+2\sqrt{t(t-x)},
\end{equation}
and
\begin{equation}\label{eqD3_7}
\varphi^{\prime\prime}(s_0)=\frac{2t^{3/2}(t-x)^{3/2}}{\sqrt{\rho}x^3}.
\end{equation}
Using (\ref{eqD3_5})-(\ref{eqD3_7}) in (\ref{eqD3_4}) leads to (\ref{eqth12}).

Next we consider $x,t\rightarrow\infty$ but with $t/x\approx 1$. The previous calculation is not valid since the saddle point $s_0\rightarrow\infty$. Now we return to the second expression in (\ref{eqD3_2}) and note that $r(s)\sim 1/s$ as $s\rightarrow\infty$. Then we approximate the integrand for $s$ large (more precisely we can scale $s=O(x)$ with $x(t-x)=O(1)$) to get
\begin{eqnarray}\label{eqD3_8}
p(t|x)&\sim& \frac{(1-\rho)e^{-\rho x}}{2\pi i}\int_{Br}e^{s(t-x)}e^{\rho x/s}ds\\
      &=& (1-\rho)e^{-\rho x}\delta(t-x)+\frac{(1-\rho)e^{-\rho x}}{2\pi i}\int_{Br}e^{s(t-x)}\left(e^{\rho x/s}-1\right)ds\nonumber\\
      &=& (1-\rho)e^{-\rho x}\delta(t-x)+(1-\rho)e^{-\rho x}\frac{d}{dt}\left[\frac{1}{2\pi i}\int_{Br}\frac{e^{s(t-x)}}{s}\left(e^{\rho x/s}-1\right)ds\right]\nonumber.
\end{eqnarray}
The last integral is equal to $I_0\left(2\sqrt{\rho x (t-x)}\right)-1$, where $I_0$ is the modified Bessel function. Differentiating this with respect to $t$ leads to (\ref{eqth11}) for $t>x$.

Now we consider $x,t\rightarrow\infty$ but with $x/t$ small. In this limit the saddle point $s_0\rightarrow s_-$, which is a branch point of the integrand in the second expression in (\ref{eqD3_2}) (but not of the first). We shift the $Br$ contour to $Br_-$, which lies slightly to the right of $s=s_-$. As $s\downarrow s_-$ we have $x(1-\rho r^2)/r\sim 2\rho^{1/4}\sqrt{s-s_-}x$, so for $x$ large we scale $s-s_*=O(x^{-2})$ by setting
\begin{equation}\label{eqD3_9}
s=-1-\rho+2\sqrt{\rho}+z/x^2=s_-+z/x^2.
\end{equation}
We then have
\begin{eqnarray}\label{eqD3_10}
x\varphi(s)&=& s(t-x)+\rho\left[\frac{x}{\sqrt{\rho}}-\rho^{-3/4}\sqrt{z}+O(x^{-1})\right]\\
           &=& -(1-\sqrt{\rho})^2t+(1+\rho-\sqrt{\rho})x-\rho^{1/4}\sqrt{z}+\frac{t}{x^2}z+O(x^{-1}).\nonumber
\end{eqnarray}
If we assume that $x=O(\sqrt{t})$ with $a=t/x^2$, then (\ref{eqD3_2}) becomes
\begin{eqnarray}\label{eqD3_11}
p(t|x)\sim && 2\rho^{-1/4}(1-\rho)e^{-(1-\sqrt{\rho})^2t}e^{(1-\sqrt{\rho})x}\\
      &\cdot& \frac{1}{2\pi i}\int_{Br_-}\frac{\sqrt{z}e^{az}e^{-\rho^{1/4}\sqrt{z}}}{(1-\sqrt{\rho})^2\left[1-e^{-2\rho^{1/4}\sqrt{z}}\right]x^3}dz.\nonumber
\end{eqnarray}
Expanding $\left(1-e^{-2\rho^{1/4}\sqrt{z}}\right)^{-1}$ as a geometric series and using the identity
\begin{eqnarray}\label{eqD3_12}
&&\frac{1}{2\pi i}\int_{Br_-}e^{az}e^{-(2n+1)\rho^{1/4}\sqrt{z}}\sqrt{z}dz\\
    &&=\frac{1}{4\sqrt{\pi}}\exp\left[-\frac{(2n+1)^2\sqrt{\rho}}{4a}\right]\left[(2n+1)^2\sqrt{\rho}-2a\right]a^{-5/2}\nonumber
\end{eqnarray}
we obtain (\ref{eqth136}). Then (\ref{eqth137}) follows upon using the Poisson summation formula in the form
\begin{equation}\label{eqD3_13}
\sum_{n=-\infty}^\infty\psi(n)=\sum_{m=-\infty}^\infty\hat{\Psi}(2\pi m)=\sum_{m=-\infty}^\infty\int_{-\infty}^\infty e^{2\pi i ym}\psi(y)dy,
\end{equation}
where $\hat{\Psi}$ is the Fourier transform of $\psi$, and $\psi(n)$ is given by the right side of (\ref{eqD3_12}). Alternately, we can set $z=w^2$ in (\ref{eqD3_11}) and after some contour deformation write the integral as
\begin{equation}\label{eqD3_14}
\frac{1}{\pi i (1-\sqrt{\rho})^2x^3}\int_{Br_*}\frac{w^2e^{aw^2}e^{-\rho^{1/4}w}}{1-e^{-2\rho^{1/4}w}}dw,
\end{equation}
where $\Re(w)>0$ on $Br_*$. The integrand in (\ref{eqD3_14}) has simple poles at $w=\rho^{-1/4}N\pi i$ for $N=\pm1,\pm2,...,$ and is anti-symmetric under the transformation $w\rightarrow-w$. We can thus write it as an infinite residue series and this again leads to (\ref{eqth137}).

Finally, we consider the case $x=O(1)$ and $t\rightarrow\infty$. We use the conformal map $s=-1-\rho+2\sqrt{\rho}\cosh\eta$. Then in the $\eta$-plane, (\ref{eqD3_2}) becomes
\begin{eqnarray}\label{eqD3_15}
p(t|x)=&&\frac{2\sqrt{\rho}(1-\rho)}{2\pi i}\\
&&\cdot\int_C \frac{e^{(-1-\rho+2\sqrt{\rho}\cosh\eta)t}e^{(1-2\sqrt{\rho}\cosh\eta+\sqrt{\rho}e^{-\eta})x}(1-e^{-2\eta})\sinh\eta}
{(1-\sqrt{\rho}e^{-\eta})^2-(\sqrt{\rho}-e^{-\eta})^2e^{-2\sqrt{\rho}x\sinh\eta}}d\eta\nonumber.
\end{eqnarray}
Here the contour $C$ goes from $\eta=\infty-i\pi/2$ to $\eta=\infty+i\pi/2$ as indicated in Figure 1. It cannot enclose any of the poles, whose location is discussed below. We let $\eta=u+iv$ where $u=\Re(\eta)$ and $v=\Im(\eta)$.

Note that $\eta=0$ is a simple zero of the integrand, since $(1-e^{-2\eta})\sinh\eta=O(\eta^2)$ and the denominator in (\ref{eqD3_15}) has a simple zero at $\eta=0$. From numerical calculations we find that the dominant singularity (which maximizes $-1-\rho+\sqrt{\rho}\cosh\eta$) occurs either along the imaginary axis for $|v|<\pi$, or along the lines $\Im(\eta)=v=\pm\pi$. In Figure 2 we sketch all the singularities of the integrand in (\ref{eqD3_15}) for $x=1$, $\rho=0.3$ (so that $x>x_*\equiv\frac{1-\sqrt{\rho}}{\sqrt{\rho}(1+\sqrt{\rho})}\doteq0.5335$). The dominant singularity is now on the imaginary axis at $\eta=\pm i \widetilde{v}$ where $\widetilde{v}\doteq 1.6111$. In Figure 2, the singularities are the crossing points of two curves, corresponding to the real and imaginary parts of the denominator in (\ref{eqD3_15}), with $\eta=u+iv$. As $x$ decreases through the value $x_*$, the dominant singularity ``turns the corner" and begins to lie on $v=\pm\pi$. In Figure 3 we sketch all the singularities of the integrand in (\ref{eqD3_15}) for $x=0.5$, $\rho=0.3$ (so that $x<x_*\doteq0.5335$), with the numerical value of the dominant singularity being $\eta=\pm i\pi+\widetilde{u}$ where $\widetilde{u}\doteq 0.5195$. In Table 1 we give the dominant singularity (with $v>0$) for various values of $x$, and various $\rho<1$. When $x=x_*$, the dominant singularity is at $\eta=\pm i\pi$. Note that when $x=x_*$, $\eta=i\pi$ is a double zero of the numerator and a triple zero of the denominator of (\ref{eqD3_15}), so $\eta=i\pi$ is a simple pole of the integrand. The Table also indicates the behavior of the singularity for $x\rightarrow 0$ and $x\rightarrow \infty$, where (\ref{eqth1410}) and (\ref{eqth1411}) apply.

Expression (\ref{eqth148}) follows by indenting the contour in Figure 1 around the two dominant singularities and using the symmetry $(v\rightarrow-v)$ of the integrand in (\ref{eqD3_15}). To obtain (\ref{eqth1411}) we use (\ref{eqth1410}) and note that $v\rightarrow 0$ as $x\rightarrow\infty$ with $x\sqrt{\rho}\sin v\rightarrow\pi$, so that $v\sim \pi/(\sqrt{\rho}x)$, which can be refined to the estimate
\begin{equation}\label{eqD3_16}
v=\frac{\pi}{\sqrt{\rho}x}-\frac{\pi(1+\sqrt{\rho})}{\rho(1-\sqrt{\rho})x^2}+O(x^{-3}).
\end{equation}
To obtain (\ref{eqth1411}) we use
$$s=-1-\rho+2\sqrt{\rho}\cos v=-(1-\sqrt{\rho})^2-\sqrt{\rho}v^2+O(v^4)$$
and (\ref{eqD3_16}).

To get the dominant singularity as $x\rightarrow 0^+$ we set $\eta=u+i\pi$ to get
\begin{equation}\label{eqD3_17}
\left(\frac{1+\sqrt{\rho}e^{-u}}{\sqrt{\rho}+e^{-u}}\right)^2=e^{2\sqrt{\rho}x\sinh u}.
\end{equation}
As $x\rightarrow 0^+$ we must have $u\rightarrow\infty$ and, using $2\sinh u\sim e^u$, we obtain from (\ref{eqD3_17}) 
\begin{equation}\label{eqD3_18}
\sqrt{\rho}xe^u=\log\left[1/\rho+O(e^{-u})\right]
\end{equation}
so that
$$u\sim -\log x+\log\log(\rho^{-1})-\frac{1}{2}\log\rho.$$
By considering higher order terms in (\ref{eqD3_18}) this can be refined to 
\begin{eqnarray*}
u= && -\log x+\log\log(\rho^{-1})-\frac{1}{2}\log\rho-\frac{2(1-\rho)}{(\log \rho)^2}x\\
      && +\frac{\rho(\log\rho)^2+(\rho^2-1)\log\rho-6(1-\rho)^2}{(\log\rho)^4}x^2+O(x^3).
\end{eqnarray*}
Then (\ref{eqth1412}) follows by using $r_*(x)=-1-\rho-\sqrt{\rho}(e^u+e^{-u})$.

To compute the unconditional density $p(t)$ we use (\ref{eqI2}) with $\mu=1$. The major contribution will come from the asymptotic matching region between the scales $x=O(1)$ and $x=O(\sqrt{t})$. In this range we use (\ref{eqth148}) with $r_*(x)$ approximated by (\ref{eqth1411}) and
\begin{equation}\label{eqD3_19}
F(x)\sim \frac{2\pi^2(1+\sqrt{\rho})}{\rho(1-\sqrt{\rho})x^3}e^{(1-\sqrt{\rho})x},  x\rightarrow\infty,
\end{equation}
which follows from the definition of $F$ below (\ref{eqth148}) and from (\ref{eqD3_16}). Scaling $x=t^{1/3}y=O(t^{1/3})$, (\ref{eqI2}) becomes asymptotically
\begin{eqnarray}\label{eqD3_20}
&&e^{-(1-\sqrt{\rho})^2t}\frac{(1+\sqrt{\rho})\pi^2}{(1-\sqrt{\rho})\rho}t^{-2/3}\\
&&\cdot\int_{x_0t^{-1/3}}^{t^{2/3}} \frac{1}{y^3}\exp\left[\frac{2\pi^2(1+\sqrt{\rho})}{\rho(1-\sqrt{\rho})y^3}\right]   \exp\left[t^{1/3}\left(-\sqrt{\rho}y-\frac{\pi^2}{\sqrt{\rho}y^2}\right)\right]dy.\nonumber
\end{eqnarray}
Here $x_0>0$ so as to avoid integrating through $y=0$, and $x_0\gg1$.

Expanding (\ref{eqD3_20}) by the Laplace method, with the major contribution coming from $y=(2\pi^2/\rho)^{1/3}$, leads to
\begin{equation}\label{eqD3_21}
p(t)\sim e^{-At}e^{-Bt^{1/3}}C_*t^{-5/6}
\end{equation}
where
$$A=(1-\sqrt{\rho})^2,$$
$$B=3\left(\frac{\pi}{2}\right)^{2/3}\rho^{1/6},$$
$$C_*=2^{2/3}3^{-1/2}\pi^{5/6}\rho^{-5/12}\frac{1+\sqrt{\rho}}{1-\sqrt{\rho}}\exp\left(\frac{1+\sqrt{\rho}}{1-\sqrt{\rho}}\right).$$
By using (\ref{eqD3_21}) and (\ref{eqI4}) we get
$$p(t)=-\frac{d}{dt}\Pr[\mathbf{V}>t]\sim \rho^{-1}AC_*t^{-5/6}e^{-At}e^{-Bt^{1/3}}$$
and this agrees with the result of Flatto (see Theorem 1.2 in \cite{FL}), if we allow for the different scaling of time used therein (where the arrival rate was taken as unity).

\section{Brief derivations for the case $\rho\approx 1$}

Now we consider the case in which the traffic intensity is close to one, and let $\rho=1-\epsilon$ with $0<\epsilon\ll 1$.

First, we consider $x=O(1)$ and $t=O(1)$. From (\ref{eqS2}) we have
$$r=r(s;1-\epsilon)=\frac{1}{2}\left[s+2-\sqrt{s^2+4s}\right]+O(\epsilon)=r_0(s)+O(\epsilon).$$
This leads to (\ref{eqth21}). On this scale the solution does not simplify much, but there is little probability mass in heavy traffic on the time scale $t=O(1)$.

Next, we consider $x=O(1)$ but large time scales $t=T/\epsilon=O(\epsilon^{-1})$. In (\ref{eqS3}), we replace $\rho$ as $1-\epsilon$ and scale $s$ as $\epsilon w$, and we have
$$ p(t|x)=\frac{\epsilon}{2\pi i}\int_{Br}e^{wT}\left[\frac{1}{1+wx}-\frac{(3+x)x^2w^2}{6(1+wx)^2}\epsilon+O(\epsilon^2)\right]dw.$$
Then (\ref{eqth22}) follows upon using the following identities
$$ \frac{1}{2\pi i}\int_{Br}\frac{e^{wT}}{1+wx}dw=\frac{1}{x}e^{-T/x},$$

$$ \frac{1}{2\pi i}\int_{Br}e^{wT}\frac{w^2}{(1+wx)^2}dw=\frac{\delta(T)}{x^2}-\frac{2x-T}{x^4}e^{-T/x}.$$
Note that by removing the condition on $x$ and expressing the result in terms of the modified Bessel function, we can recover the result of Morrison (see (3.12) in \cite{MO}). The term proportional to $\delta(T)$ does not actually mean there is mass at $T=0$, but rather corresponds to the small ($O(\epsilon)$) mass that exists in the shorter time scale $t$.

Now consider $x=X/\epsilon=O(\epsilon^{-1})$ and $t=T/\epsilon=O(\epsilon^{-1})$ with $1<T/X<\infty$. We define
\begin{eqnarray*}
\phi(s)&=& s(\frac{T}{X}-1)+(1-\epsilon)\left(r(s;1-\epsilon)-1\right)\\
      &=& s\frac{T}{X}-\frac{1}{2}\left(s+\sqrt{s^2+4s}\right)+\frac{1}{2}\left(1+\frac{s}{\sqrt{s^2+4s}}\right)\epsilon+O(\epsilon^2)\\
      &=&\phi_0(s)+\phi_1(s)\epsilon+O(\epsilon^2)
\end{eqnarray*}
and rewrite (\ref{eqS3}) as
\begin{eqnarray}\label{eqD4_1}
p(t|x)&=& \frac{\epsilon}{2\pi i}\int_{Br}\frac{e^{\phi(s)X/\epsilon}\left[1-(1-\epsilon)r^2\right]}{\left[1-(1-\epsilon)r\right]^2-(1-\epsilon)(1-r)^2e^{-x\left[1-(1-\epsilon)r^2\right]/r}}ds\nonumber\\
      &\sim& \frac{\epsilon}{2\pi i}\int_{Br}\frac{1-(1-\epsilon)r^2}{\left[1-(1-\epsilon)r\right]^2}e^{\phi(s)X/\epsilon}
ds\nonumber\\
      &\sim& \frac{\epsilon}{2\pi i}\int_{Br}\sqrt{1+\frac{4}{s}}e^{\phi_1(s)X}e^{\phi_0(s)X/\epsilon}ds.
\end{eqnarray}
The saddle point in the last approximation in (\ref{eqD4_1}), where $\phi_0^\prime(s)=0$, is $s=s_0^*$ where 
\begin{equation}\label{eqD4_2}
s_0^*=s_0^*(\frac{T}{X})=-2+\frac{2T-X}{\sqrt{T(T-X)}}.
\end{equation}
Note that in the heavy traffic case, $s_-=-(1-\sqrt{\rho})^2=O(\epsilon^2)$, $s_-<s_0^*<\infty$ whenever $1<T/X<\infty$, $s_0^*\rightarrow\infty$ as $T/X\downarrow 1$, while $s_0^*\rightarrow 0$ as $T/X\rightarrow\infty$. By shifting $Br$ in (\ref{eqD4_1}) to $Br^{\prime\prime}$, on which $\Re(s)=s_0^*$, and using the standard Laplace method, we have
$$ p(t|x)\sim\frac{\epsilon^{3/2}\sqrt{1+4/s_0^*}}{\sqrt{2\pi}\sqrt{X\phi_0^{\prime\prime}(s_0^*)}}e^{\phi_1(s_0^*)X}e^{\phi_0(s_0^*)X/\epsilon}. $$
This leads to (\ref{eqth24}). We note that if $X\rightarrow\infty$ and $T\rightarrow\infty$ but $T/X=O(1)$, the approximation (\ref{eqth24}) remains valid.

For the case $X=O(1)$, $T=O(1)$ and $T-X\rightarrow 0^+$, we notice that the saddle point $s_0^*\rightarrow\infty$. Then by the same argument as in the case when $\rho$ is fixed and less than one, we can easily get (\ref{eqth23}), where $T_*=(T-X)/\epsilon^2=O(1)$.

Next, we consider $X=\sqrt{\epsilon}Z=O(\sqrt{\epsilon})$, $T=O(1)$. It follows that $T/Z^2=O(1)$ and the saddle point $s_0^*\rightarrow 0$. We shift the contour $Br$ to $Br_0$, which lies slightly to the right of $s=0$. By scaling $s=O(\epsilon)=\epsilon w$ in the first equation in (\ref{eqD4_1}), we have
$$ p(t|x)\sim \frac{2\epsilon^{3/2}}{2\pi i}\int_{Br_0}\frac{e^{wT}e^{-\sqrt{w}Z}}{\sqrt{w}\left(1-e^{-2\sqrt{w}Z}\right)}dw. $$
Expanding $\left(1-e^{-2\sqrt{w}Z}\right)^{-1}$ as a geometric series and using the identity
$$ \frac{1}{2\pi i}\int_{Br_0}\frac{e^{wT}e^{-(2n+1)\sqrt{w}Z}}{\sqrt{w}}dw=\frac{1}{\sqrt{\pi T}}e^{-(2n+1)^2Z^2/(4T)},$$
we obtain (\ref{eqth25}).
Then (\ref{eqth25psf}) follows by using the Poisson summation formula (\ref{eqD3_13}). Note that from (\ref{eqth25psf}), we can easily verify that cases 5 and 2 in Theorem 2 asymptotically match, in the intermediate limit where $x\rightarrow\infty$ and $Z\rightarrow 0$. Similarly, cases 4 and 5 match in the limit $X\rightarrow 0$ and $Z\rightarrow\infty$, which follows easily from (\ref{eqth25}).

Now we consider $X=O(1)$ and $T=\Theta/\epsilon=O(\epsilon^{-1})$ (thus $x=O(\epsilon^{-1})$ and $t=O(\epsilon^{-2})$). In the first expression in (\ref{eqD4_1}), we scale $s=O(\epsilon^{2})$ by setting $s=(\xi-1/4)\epsilon^2$, which leads to (\ref{eqth26a}) after some calculation.

Furthermore, we expand the integrand in (\ref{eqth26a}) as a geometric series, and we have 
\begin{equation}\label{eqD4_3}
p(t|x)\sim \frac{4\epsilon^2e^{X/2}e^{-\Theta/4}}{\pi i}\int_{Br}\frac{e^{\xi\Theta}\sqrt{\xi}}{(1+2\sqrt{\xi})^2}\sum_{n=0}^\infty \left(\frac{1-2\sqrt{\xi}}{1+2\sqrt{\xi}}\right)^{2n}e^{-(2n+1)X\sqrt{\xi}}d\xi.
\end{equation}
Note that if we let $X\rightarrow\infty$, then the $n=0$ term in (\ref{eqD4_3}) dominates, and we have
\begin{eqnarray*}
&&\frac{1}{2\pi i}\int_{Br}e^{\xi\Theta}\frac{\sqrt{\xi}}{(1+2\sqrt{\xi})^2}e^{-X\sqrt{\xi}}d\xi\\
&=&\frac{e^{-X^2/(4\Theta)}}{4\sqrt{\pi\Theta}}+\frac{\sqrt{\Theta}}{8\sqrt{\pi}}e^{-X^2/(4\Theta)}-\frac{4+X+\Theta}{8}e^{\Theta/4}e^X\mathrm{erfc}\left(\frac{X+\Theta}{2\sqrt{\Theta}}\right)\\
&\sim& \frac{e^{-X^2/(4\Theta)}}{4\sqrt{\pi\Theta}},\textrm{ }\textrm{ }\textrm{ } X\rightarrow\infty.
\end{eqnarray*}
Here we used 
$$\mathrm{erfc}(z)=1-\mathrm{erf}(z)=\frac{2}{\sqrt{\pi}}\int_z^\infty e^{-t^2}dt\sim\frac{1}{\sqrt{\pi}z}e^{-z^2},\textrm{ as }z\rightarrow\infty.$$
Then (\ref{eqD4_3}) becomes, for $\Theta$ fixed and $X\rightarrow\infty$,
\begin{equation}\label{eqD4_4}
p(t|x)\sim\frac{2\epsilon^2}{\sqrt{\pi\Theta}}e^{X/2}e^{-\Theta/4}e^{-X^2/(4\Theta)}.
\end{equation}
When $X=O(\epsilon^{-1})$ and $T=O(\epsilon^{-1})$ but with $T/X=O(1)$, (\ref{eqth24}) remains valid, and letting $T/X\rightarrow\infty$ in (\ref{eqth24}) regains (\ref{eqD4_4}). This again verifies that these two cases asymptotically match.

We return to (\ref{eqD4_3}), and let $\sqrt{\xi}=(z-1)/2$, with which the integral becomes
\begin{equation}\label{eqD4_5}
\frac{e^{\Theta/4}}{2}\int_{C^\prime}(z-1)^2\sum_{n=0}^\infty\frac{(z-2)^{2n}}{z^{2n+2}}e^{(2n+1)X/2}e^{-[(2n+1)X+\Theta]z/2}e^{\Theta z^2/4}dz.
\end{equation}
Here the contour $C^\prime$ can be taken as the imaginary axis in the $z$-plane, indented to the right of $z=0$. Using the binomial expansion on $(z-2)^{2n}$
$$(z-2)^{2n}=\sum_{m=0}^{2n}\frac{(-1)^m(2n)!}{m!(2n-m)!}2^{2n-m}z^m $$
and after some simplification, (\ref{eqD4_3}) leads to 
\begin{eqnarray}\label{eqD4_6}
p(t|x)\sim&& 2\epsilon^2\sum_{n=0}^\infty\sum_{m=0}^{2n}\frac{(-1)^m(2n)!}{m!(2n-m)!}2^{2n-m}e^{(n+1)X}\nonumber\\
&\cdot& \frac{1}{2\pi i}\int_{C^\prime}(z-1)^2\frac{e^{-A_nz}e^{\Theta z^2/4}}{z^N}dz,
\end{eqnarray}
where
$$ A_n=\frac{2n+1}{2}X+\frac{\Theta}{2}$$
and
$$ N=2n-m+2.$$
We express the integral in (\ref{eqD4_6}) in terms of parabolic cylinder functions, using
$$ \frac{1}{2\pi i}\int_{Br}z^\nu e^{-wz+z^2/2}dz=\frac{1}{\sqrt{2\pi}}D_\nu(w)e^{-w^2/4}, $$
thus obtaining
\begin{eqnarray}\label{eqD4_7}
&&p(t|x)\sim \frac{\epsilon^2}{\sqrt{2\pi}} \sum_{n=0}^\infty e^{(n+1)X}e^{-z_n^2/4}\sum_{m=0}^{2n}\frac{(-1)^m(2n)!}{m!(2n-m)!}(2\Theta)^{(2n-m)/2}\nonumber\\
&&\cdot\left[\frac{4}{\sqrt{2\Theta}}D_{m-2n}(z_n)-4D_{m-2n-1}(z_n)+\sqrt{2\Theta}D_{m-2n-2}(z_n)\right],
\end{eqnarray}
where $$z_n=\frac{A_n}{\sqrt{\Theta/2}}=\frac{(2n+1)X+\Theta}{\sqrt{2\Theta}}.$$
Replacing $2n-m$ by $l$, (\ref{eqD4_7}) leads to (\ref{eqth26b}).

If we let $X\rightarrow 0$ and $\Theta\rightarrow 0$ with $X/\sqrt{\Theta}$ (thus $z_n$) fixed, the term with $l=0$ in (\ref{eqth26b}) dominates and we have
\begin{eqnarray}\label{eqD4_8}
p(t|x) &\sim& \frac{2\epsilon^2}{\sqrt{\pi\Theta}}\sum_{n=0}^\infty e^{(n+1)X}e^{-z_n^2/4}D_0(z_n)\nonumber\\
       &\sim& \frac{2\epsilon^2}{\sqrt{\pi\Theta}}\sum_{n=0}^\infty e^{-(2n+1)^2X^2/(4\Theta)}.
\end{eqnarray}
Here we note that $X^2/\Theta=Z^2/T$ and we used the fact that $D_0(w)=e^{-w^2/4}$. Since (\ref{eqD4_8}) is the same as (\ref{eqth25}), we have shown that case 5 is really a special case of case 6 in Theorem 2. 

Alternately, we can treat the problem on the $(X,\Theta)$ scale by using the conformal map in section 2 to get (\ref{eqD3_15}). The poles of the integrand in (\ref{eqD3_15}) satisfy
$$\left(\frac{1-\sqrt{\rho}e^{-\eta}}{\sqrt{\rho}-e^{-\eta}}\right)^2=\exp\left(\frac{-2\sqrt{\rho}X\sinh\eta}{\epsilon}\right).$$
Since $\rho=1-\epsilon$, we must scale $\eta$ to be $O(\epsilon)$. Let $\eta=\eta_1\epsilon+\eta_2\epsilon^2+O(\epsilon^3)$, we find that 
\begin{equation}\label{eqD4_9}
\left(\frac{2\eta_1+1}{2\eta_1-1}\right)^2=e^{-2X\eta_1}
\end{equation}
and $\eta_2=\eta_1/2$.
Setting $\eta_1=u+iv$ we find that all the roots of (\ref{eqD4_9}) are on the imaginary axis, and with $u=0$ (\ref{eqD4_9}) becomes 
$$ (4v^2-1)\sin(Xv)=4v\cos(Xv) $$
or $$ \left[\cos\left(\frac{Xv}{2}\right)-2v\sin\left(\frac{Xv}{2}\right)\right]\left[2v\cos\left(\frac{Xv}{2}\right)+\sin\left(\frac{Xv}{2}\right)\right]=0.$$
Thus we have either $\cot(Xv/2)=2v$ or $\tan(Xv/2)=-2v$. Denoting the $n^\mathrm{th}$ positive solution by $v_n=v_n(X)$ we have $v_n(X)\in\left((n-1)\pi/X,n\pi/X\right)$, $n\geq 1$. The first solution satisfies $\cot(Xv_1/2)=2v_1$ and has the asymptotic expansions
$$ v_1(X)=\frac{\pi}{X}-\frac{4\pi}{X^2}+\frac{16\pi}{X^3}+O(X^{-4}), \textrm{ }X\rightarrow\infty,$$
$$ v_1(X)=\frac{1}{\sqrt{X}}-\frac{\sqrt{X}}{24}+\frac{11}{5760}X^{3/2}+O(X^{5/2}), \textrm{ }X\rightarrow 0.$$
The higher order solutions behave as
$$ v_n(X)=n\pi\left[\frac{1}{X}-\frac{4}{X^2}+\frac{16}{X^3}\right]+O(X^{-4}), \textrm{ }X\rightarrow\infty,$$
and
$$ v_n(X)=(n-1)\frac{\pi}{X}+\frac{1}{(n-1)\pi}-\frac{X}{(n-1)^3\pi^3}+O(X^2), \textrm{ }n\geq 2, X\rightarrow 0.$$
Then in the $s$-plane the singularities satisfy $s=O(\epsilon^2)$ with
\begin{eqnarray}\label{eqD4_10}
\widetilde{r}_*(X;v_n)&=& (-1-\rho+2\sqrt{\rho}\cosh\eta)/\epsilon^2\nonumber\\
                  &=& -(v_n^2+1/4)-(v_n^2+1/4)\epsilon/2+O(\epsilon^2),
\end{eqnarray}
and we also have, at the singular points,
\begin{equation}\label{eqD4_11}
\exp\left[\left(1-2\sqrt{\rho}\cosh\eta+\sqrt{\rho}e^{-\eta}\right)X/\epsilon\right]\sim e^{X/2-iXv_n}.
\end{equation}
Expression (\ref{eqth26c}) follows by indenting the integration contour in (\ref{eqD3_15}) around all the singularities, using the symmetry ($v\rightarrow -v$) of the integrand and also using (\ref{eqD4_10}) and (\ref{eqD4_11}). Note that in heavy traffic all the singularities are close to $\eta=0$ and on the imaginary axis. If we retain only the leading term in $\widetilde{r}_*$, (\ref{eqth26c}) is equivalent to (\ref{eqth26a}), as can be seem by setting $\xi=\eta_1^2$ in (\ref{eqth26a}), with which the singularities satisfy (\ref{eqD4_9}).

If we consider even larger time scales, with $t=\sigma/\epsilon^3=O(\epsilon^{-3})$ (thus $t=O(\epsilon^{-3})$ ), then the two singularities at $\eta_1=\pm iv_1$ dominate. Here $v_1=v_1(X)$ is the unique root in the interval $(0,\pi/X)$ of (\ref{eqth26croot}). This leads to (\ref{eqth26cmatch}). Note that to compute the leading term on the $\sigma$-time scale, we need the $O(\epsilon)$ correction term to $\widetilde{r}_*$.

To compute the unconditional density $p(t)$ on the scale $t=O(\epsilon^{-3})=\sigma/\epsilon^3$, we use (\ref{eqI2}) with $\mu=1$ and (\ref{eqth26cmatch}). For simplification, we write $V$ for $v_1(X)$, and we have 
\begin{eqnarray}\label{eqD4_12}
p(t)&\sim& \epsilon\int_0^{\epsilon t}\widetilde{F}(X;V)e^{\widetilde{r}_*(X;V)\sigma/\epsilon}e^{-X/\epsilon}dX\nonumber\\
      &\sim& \epsilon\int_0^\infty \widetilde{F}(X;V)e^{-(V^2+1/4)\sigma/2}e^{\left[-X-(V^2+1/4)\sigma\right]/\epsilon}dX.
\end{eqnarray}
This is a Laplace type integral. We let $\Phi(X)=-X-(V^2+1/4)\sigma$ and must find the maximum of this function over $X$. Then $\Phi^\prime(X)=-1-2\sigma V(X)V^\prime(X)$, and from (\ref{eqth26croot}) by implicit differentiation we obtain 
$$ V^\prime(X)=-\frac{(1+4V^2)V}{(1+4V^2)X+4}.$$
The condition $\Phi^\prime(X)=0$ defines $X=X(\sigma)$ implicitly, where the inverse function is 
\begin{equation}\label{eqD4_13}
\sigma=\frac{X(1+4V^2)+4}{2V^2(1+4V^2)},
\end{equation}
where $V=V(X)$. Denoting the right-hand side of (\ref{eqD4_13}) as $k(X)$, we can verify that $k^\prime(X)>0$, so that $k(X)$ is a monotonically increasing function, and that $k(X)\rightarrow\infty$ as $X\rightarrow\infty$. From our asymptotic results $V\sim1/\sqrt{X}$ as $X\rightarrow 0$, so that $k(X)\rightarrow 0$ as $X\rightarrow 0$. Hence there is a unique positive root of $k(X)=\sigma$, which we denote by $\widetilde{X}$, so that 
\begin{equation}\label{eqD4_14}
\sigma=\frac{\widetilde{X}(1+4\widetilde{V}^2)+4}{2\widetilde{V}^2(1+4\widetilde{V}^2)}, \textrm{ }\widetilde{V}=V(\widetilde{X}).
\end{equation}
Then we use the standard Laplace method in (\ref{eqD4_12}) and get
\begin{equation}\label{eqD4_15}
p(t)\sim\frac{\sqrt{2\pi}\epsilon^{3/2}}{\sqrt{-\Phi^{\prime\prime}(\widetilde{X})}}\widetilde{F}(\widetilde{X};\widetilde{V})e^{-(\widetilde{V}^2+1/4)\sigma/2}e^{\Phi(\widetilde{X})/\epsilon}.
\end{equation}

The first root $v_1(X)=V$ satisfies
\begin{equation}\label{eqD4_16}
\cot\left(\frac{XV}{2}\right)=2V.
\end{equation}
Letting $\psi=\widetilde{X}\widetilde{V}$, we write (\ref{eqD4_16}) in parametric form, with $\widetilde{X}=2\psi\tan(\psi/2)$ and $\widetilde{V}=\frac{1}{2}\cot(\psi/2)$. Then (\ref{eqD4_14}) becomes 
\begin{equation}\label{eqD4_17}
\sigma=4\tan^3\left(\frac{\psi}{2}\right)(\psi+\sin\psi)
\end{equation}
and we can rewrite (\ref{eqD4_15}) in terms of $\sigma$ and $\psi$, using 
\begin{equation}\label{eqD4_18}
-\Phi(\widetilde{X})=-\widetilde{X}-(\widetilde{V}+1/4)\sigma=2\psi\tan\left(\frac{\psi}{2}\right)+\frac{\sigma}{4}\csc^2\left(\frac{\psi}{2}\right)\equiv F_0(\psi),
\end{equation}
\begin{equation}\label{eqD4_19}
\left(\frac{\widetilde{V}^2}{2}+\frac{1}{8}\right)\sigma-\frac{\widetilde{X}}{2}=-\psi\tan\left(\frac{\psi}{2}\right)+\frac{\sigma}{8}\csc^2\left(\frac{\psi}{2}\right)\equiv F_1(\psi),
\end{equation}
\begin{equation}\label{eqD4_20}
\widetilde{F}(\widetilde{X};\widetilde{V})e^{-\widetilde{X}/2}=\frac{\cos^3\left(\frac{\psi}{2}\right)}{\sin\left(\frac{\psi}{2}\right)(\psi+\sin\psi)},
\end{equation}
and
\begin{equation}\label{eqD4_21}
\Phi^{\prime\prime}(\widetilde{X})=-\frac{3\cot\left(\frac{\psi}{2}\right)(\psi+\sin\psi)+4\cos^4\left(\frac{\psi}{2}\right)}{2(\psi+\sin\psi)^2}.
\end{equation}
The expressions in (\ref{eqD4_17})-(\ref{eqD4_19}) are the same as those in Morrison \cite{MO} (corresponding respectively to (A.21), (A.16) and (A.17) in \cite{MO}).
Using (\ref{eqD4_17})-(\ref{eqD4_21}) in (\ref{eqD4_15}), we have
\begin{equation}\label{eqD4_22}
p(t)\sim\frac{\sqrt{2\pi}\epsilon^{3/2}\cot\left(\frac{\psi}{2}\right)}{\sqrt{F^{\prime\prime}_0(\psi)}}e^{-F_0(\psi)/\epsilon}e^{-F_1(\psi)}
\end{equation}
and this agrees with the result of Morrison (see (3.35) in \cite{MO}) if we note that $$p(t)=-\frac{d}{dt}\Pr[\mathbf{V}>t]\sim\epsilon^2 F_0^\prime(\psi) \psi^\prime(\sigma) \Pr[\mathbf{V}>t].$$

Finally, we discuss the connection between the unconditional tail expansion of Flatto (for fixed $\rho<1$, see (1.3) in \cite{FL}) and the heavy traffic result of Morrison (see (3.35) in \cite{MO}). If $\rho\approx 1$ we have
\begin{equation}\label{eqD4_23}
\Pr\left[\mathbf{V}_{\mathrm{PS}}>t\right]=\frac{1}{\rho}\Pr\left[\mathbf{W}_{\mathrm{ROS}}>t\right]\sim\Pr\left[\mathbf{W}_{\mathrm{ROS}}>t\right].
\end{equation}
We expand the $\sigma$-scale result in \cite{MO} for $\sigma\rightarrow\infty$. It follows that $\psi\rightarrow\pi$ as $\sigma\rightarrow\infty$, with
\begin{equation}\label{eqD4_24}
\sigma\sim 4\pi \cot^3\left(\frac{\pi-\psi}{2}\right)\sim\frac{32\pi}{(\pi-\psi)^3}.
\end{equation}
Thus $\psi\sim \pi-(32\pi)^{1/3}\sigma^{-1/3}$. We replace $t$ by $\mu t$ in (3.35) in \cite{MO}, so that $\sigma=\epsilon^3\mu t$. Then for $\sigma=\epsilon^3\mu t\rightarrow\infty$ the results in \cite{MO} become 
$$ F_0(\psi)\sim \frac{\epsilon^3}{4}\mu t+3\left(\frac{\pi}{2}\right)^{2/3}\epsilon(\mu t)^{1/3}-4, $$
$$ F_1(\psi)\sim \frac{\epsilon^3}{8}\mu t-\frac{1}{2}\left(\frac{\pi}{2}\right)^{2/3}\epsilon(\mu t)^{1/3}+2, $$
and
$$ 2(\sin\psi)\left(\frac{2\pi}{\epsilon F^{\prime\prime}_0(\psi)}\right)^{1/2}\sim 2^{14/3}3^{-1/2}\pi^{5/6}\epsilon^{-3}(\mu t)^{-5/6}. $$
Thus for $\sigma\rightarrow\infty$ the heavy traffic result in \cite{MO} becomes
\begin{equation}\label{eqD4_25}
\Pr[\mathbf{V}>t]\sim\frac{\alpha^*\exp\left[-\beta^*\mu t-\gamma^*(\mu t)^{1/3}\right]}{(\mu t)^{5/6}},
\end{equation}
where
$$ \alpha^*=2^{14/3}3^{-1/2}\pi^{5/6}\epsilon^{-3}e^{4/\epsilon-2}, $$
$$ \beta^*=-\left(\frac{\epsilon^2}{4}+\frac{\epsilon^3}{8}\right), $$
$$ \gamma^*=3\left(\frac{\pi}{2}\right)^{2/3}\left(1-\frac{\epsilon}{6}\right). $$

If we let $\rho=1-\epsilon$ and replace $t$ according to $\lambda t=(1-\epsilon)\mu t$ in Flatto's result ((1.3) in \cite{FL}), and expand that expansion for $\epsilon\rightarrow 0$, we also get (\ref{eqD4_25}). This shows that the interesting tail structure in (1.3) is also inherent in the results in \cite{MO}, if we consider the latter for times even greater than the largest natural time scale.

\section{Discussion}
To summarize, we have derived several asymptotic formulas for the conditional sojourn time distribution in the $M/M/1$-PS model. Our analysis considers all of the space/time scales inherent to the problem, and leads to qualitative insights into the structure of the distribution. For example, if $\rho<1$, $x$ is large and $t/x\sim 1/(1-\rho)$, formula (\ref{eqth12}) simplifies to the Gaussian 
\begin{equation}\label{eq5_1}
p(t|x)\approx \frac{(1-\rho)^{3/2}}{2\sqrt{\rho\pi x}}\exp\bigg[-\frac{(1-\rho)^3}{4\rho x}\left(t-\frac{x}{1-\rho}\right)^2\bigg],
\end{equation}
which gives the spread about the well known mean value $\mathbf{E}\left[\mathbf{V}(x)\right]=x/(1-\rho)$.

In contrast, in the heavy traffic case, for $x=O(1)$ (short to moderate jobs) the sojourn time distribution is approximately exponentially distributed on the large time scale $T$, as shown by (\ref{eqth22}). We can show that this exponential distribution still holds for $x=Z/\sqrt{\epsilon}=O(\epsilon^{-1/2})$ (where now the sojourn time must be scaled to be $O(\epsilon^{-3/2})$), but for very large jobs, where $x=X/\epsilon=O(\epsilon^{-1})$, the distribution is concentrated on time scales $t=\Theta/\epsilon^2=O(\epsilon^{-2})$, where part 6 of Theorem 2 applies. We also note that when the condition on $x$ is removed, which corresponds to multiplying $e^{-T/x}$ by $e^{-x}$ and integrating, we obtain the modified Bessel functions that appeared in \cite{MO}.

Our results easily reproduce the unconditional asymptotics in \cite{MO} and \cite{FL}, and lead to a better understanding of what ranges of $x$ are important for removing the condition via (\ref{eqI2}). We are presently investigating the more general $M/G/1$-PS model, and our preliminary results suggest that the basic scales for the conditional sojourn time asymptotics are very similar to those here, but that the unconditional asymptotics are quite sensitive to the tail behavior of the service distribution.

\newpage
\begin{center}
\begin{table}
\caption{The dominant singularity ($\eta=u+iv$ with $v>0$) for various values of $x$ and $\rho<1$.}\label{table}
  \vspace{0.5cm}
\begin{tabular} {|c||c|c|c|c|c|}       \hline
\backslashbox{$\rho$}{$x$}     & $0.01$          & $0.05$          & $0.1$           & $0.5$           & $1.0$ \\\hline
$0.1$  & $6.5871+i\pi$ & $4.9639+i\pi$ & $4.2532+i\pi$ & $2.4866+i\pi$ & $1.5299+i\pi$\\\hline
$0.2$  & $5.8796+i\pi$ & $4.2448+i\pi$ & $3.5186+i\pi$ & $1.5625+i\pi$ & $2.3370i$\\\hline
$0.3$  & $5.3831+i\pi$ & $3.7333+i\pi$ & $2.9861+i\pi$ & $0.5195+i\pi$ & $1.6111i$\\\hline
$0.4$  & $4.9614+i\pi$ & $3.2909+i\pi$ & $2.5125+i\pi$ & $1.9978i$     & $1.2560i$\\\hline
$0.5$  & $4.5641+i\pi$ & $2.8626+i\pi$ & $2.0323+i\pi$ & $1.5312i$     & $1.0129i$\\\hline
$0.6$  & $4.1575+i\pi$ & $2.4049+i\pi$ & $1.4685+i\pi$ & $1.2118i$     & $0.8218i$\\\hline
$0.7$  & $3.7037+i\pi$ & $1.8498+i\pi$ & $0.4973+i\pi$ & $0.9534i$     & $0.6566i$\\\hline
$0.8$  & $3.1312+i\pi$ & $0.9559+i\pi$ & $1.8187i$     & $0.7188i$     & $0.5006i$\\\hline
$0.9$  & $2.1965+i\pi$ & $1.6812i$     & $1.1094i$     & $0.4748i$     & $0.3334i$\\\hline\hline

\backslashbox{$\rho$}{$x$}      & $2.0$     & $3.0$     & $4.0$     & $5.0$     & $10.0$\\\hline
$0.1$ & $2.2245i$ & $1.5826i$ & $1.2849i$ & $1.0982i$ & $0.6693i$\\\hline
$0.2$ & $1.3682i$ & $1.0513i$ & $0.8736i$ & $0.7553i$ & $0.4694i$\\\hline
$0.3$ & $1.0394i$ & $0.8127i$ & $0.6810i$ & $0.5920i$ & $0.3729i$\\\hline
$0.4$ & $0.8359i$ & $0.6599i$ & $0.5562i$ & $0.4855i$ & $0.3099i$\\\hline
$0.5$ & $0.6860i$ & $0.5454i$ & $0.4619i$ & $0.4049i$ & $0.2620i$\\\hline
$0.6$ & $0.5633i$ & $0.4506i$ & $0.3834i$ & $0.3373i$ & $0.2216i$\\\hline
$0.7$ & $0.4544i$ & $0.3655i$ & $0.3123i$ & $0.2759i$ & $0.1842i$\\\hline
$0.8$ & $0.3492i$ & $0.2823i$ & $0.2423i$ & $0.2149i$ & $0.1461i$\\\hline
$0.9$ & $0.2342i$ & $0.1903i$ & $0.1641i$ & $0.1461i$ & $0.1012i$\\\hline

\end{tabular}\\
\end{table}
\end{center}

\newpage
\begin{figure}[hbp]
\begin{center}
\includegraphics[angle=0, width=0.6\textwidth]{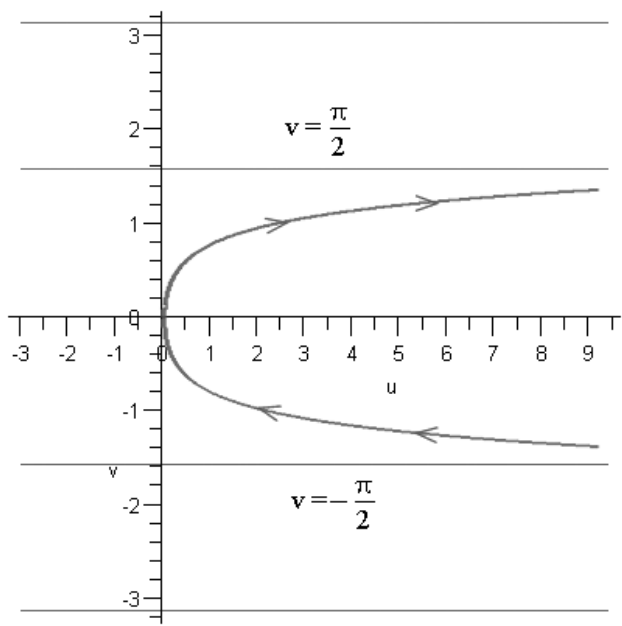}
\caption{The contour $C$.} \label{figure_c}
\end{center}
\end{figure}

\begin{figure}[hbp]
\begin{center}
\includegraphics[angle=0, width=0.8\textwidth]{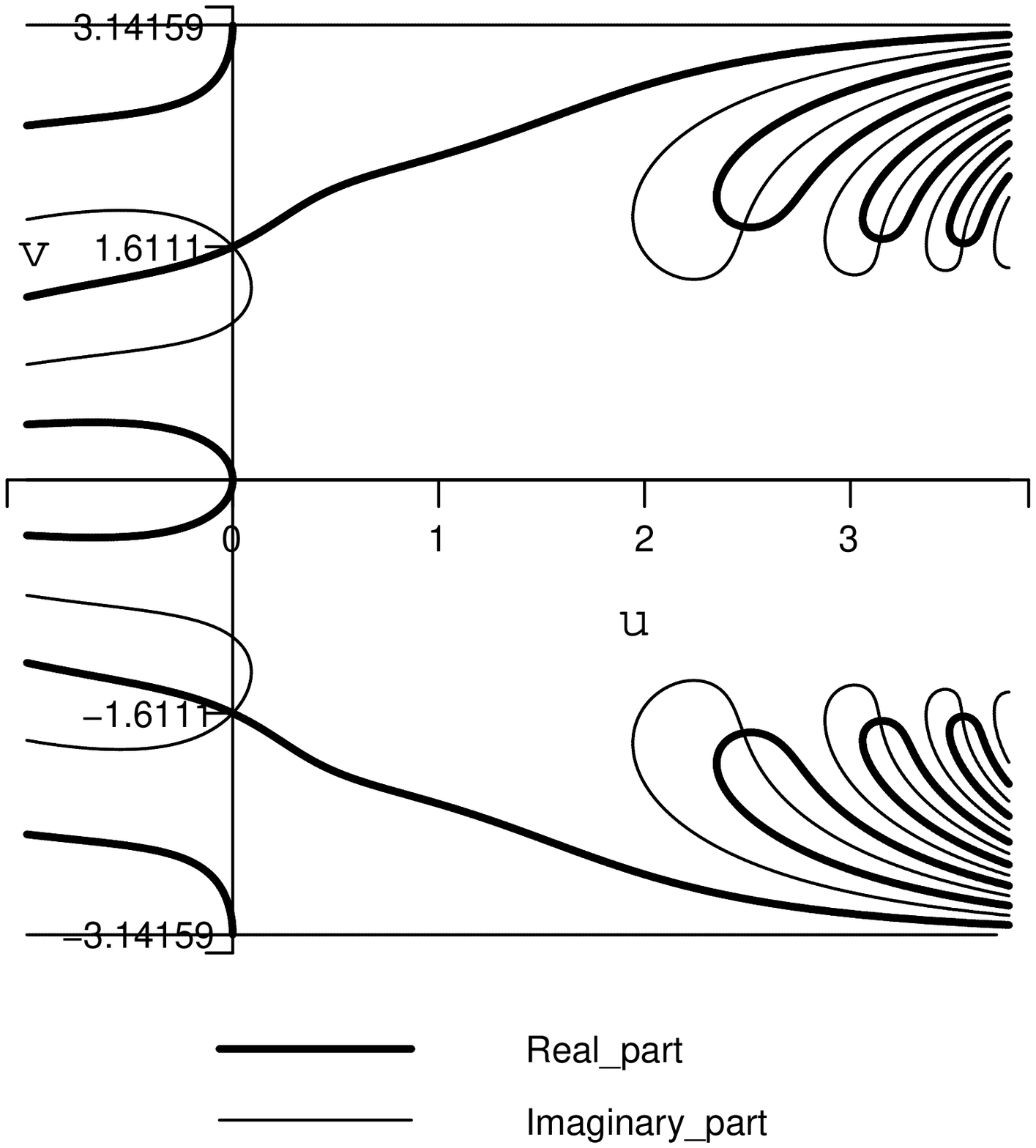}
\caption{The singularities of the integrand in (\ref{eqD3_15}) for $x=1$, $\rho=0.3$.} \label{figure_x=1}
\end{center}
\end{figure}

\begin{figure}[hbp]
\begin{center}
\includegraphics[angle=0, width=0.8\textwidth]{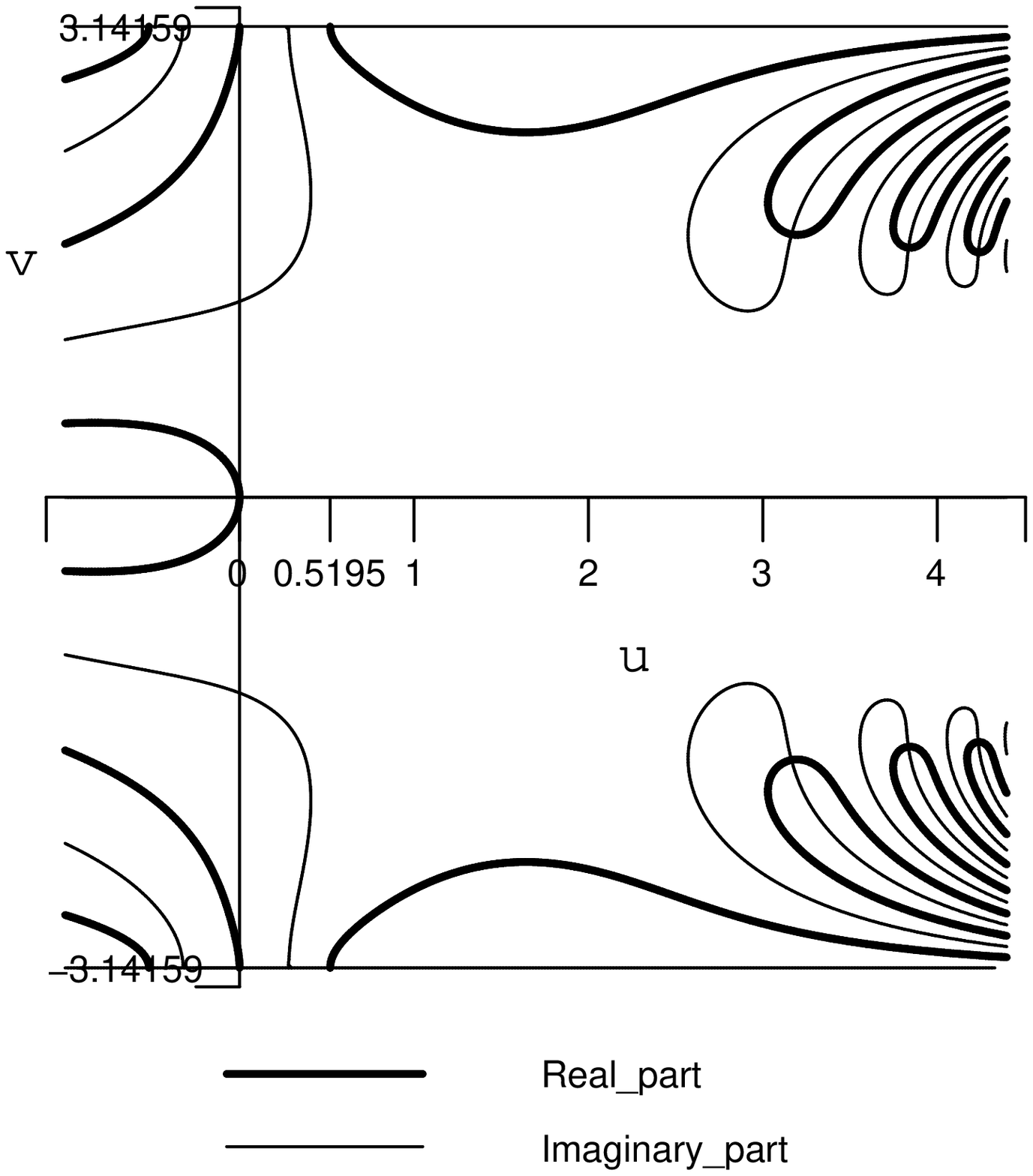}
\caption{The singularities of the integrand in (\ref{eqD3_15}) for $x=0.5$, $\rho=0.3$.} \label{figure_x=0.5}
\end{center}
\end{figure}

\end{document}